\documentclass[12pt,a4paper]{article}
\usepackage{amsthm}
\usepackage{amssymb} 
\usepackage{newpxtext,newpxmath}
\begin{document}
{\Large A general expression of the triplets of integer sided triangles with a $120^{\circ}$ angle, in parallel with the case of a $60^{\circ}$ angle.}\\\\
\centerline{Yasushi Ieno (vei04156@nifty.com)}\\\\
\centerline{Abstract}

A research for a general expression of the triplets of integer sided triangles with a $120^{\circ}$ angle, in parallel with the case of a $60^{\circ}$ angle.\\

\noindent 0.Introduction\\

Recently we have got interested in integer sided triangles after reading Zimhoni's paper [1], in which Zimhoni treated the triplets of integer sided triangles with a $60^{\circ}$ angle, showed a general expression of the triplets proved by Katayama [5], and proved another  general expression by the hint of [2] and [5].

Then we have tried to find a general expression of the triplets of integer sided triangles with a $120^{\circ}$ angle, and shall show it in this paper.\\
\\
\noindent 1. For the case of  a $60^{\circ}$ angle\\

For the integer sided triangles with a $60^{\circ}$ angle, the triplets of the sides of which are known as Eisensteinian triplets, and Zimhoni show that any such triangle is obtained by pre-multiplication (7,8,5) or (13,15,7) by a product of five matrices [1], by the hint of  Berggren [2]. 

Zimhoni defined a definition and proved a theorem in [1], then we show them as below:\\

\noindent Definition 1. 

A primitive Eisensteinian triangle is an integer sided triangle with side lengths (a,b,c) such that the angle opposite the side of length a is of $60^{\circ}$ and gcd(a,b,c){=}1. Here, we will also follow the convention of ordering the triplet such that b${>}$c.\\

\noindent Theorem 1. 

An integer sided triangle is primitive Eisensteinian if, and only if it has side lengths (1,1,1) or it has side lengths (a,b,c) or (a,b,b${-}$c) such that:\\
\[\rm
  \left(\rm
    \begin{array}{c}
    \rm a \\
    \rm b \\
    \rm c 
    \end{array}
  \right) = M\cdot v,\ where\ v{\in}
\left\{
\left(
    \begin{array}{c}
      7  \\
      8  \\
      5 
    \end{array}
  \right){,}
\left(
    \begin{array}{c}
      13  \\
      15  \\
       7 
    \end{array}
  \right)
  \right\}
\]

And M is a finite (or empty) product of the matrices:
\[\rm
  M_1=\left(
    \begin{array}{ccc}
      7 & -6 & 6 \\
      8 & -7 & 7 \\
      4 & -4 & 3
    \end{array}
  \right){,}
M_2=\left(
    \begin{array}{ccc}
      7 & 6 & -6 \\
      8 & 7 & -7 \\
      4 & 3 & -4
    \end{array}
  \right){,}
M_3=\left(
    \begin{array}{ccc}
      7 & 6 & 0 \\
      8 & 7 & 0 \\
      4 & 3 & 1
    \end{array}
  \right){,}
\]
\[\rm
  M_4=\left(
    \begin{array}{ccc}
      7 & 0 & 6 \\
      8 & 0 & 7 \\
      4 & 1 & 3
    \end{array}
  \right){,}
M_5=\left(
    \begin{array}{ccc}
      7 & 0 & -6 \\
      8 & 0 & -7 \\
      4 & 1 & -4
    \end{array}
  \right)
\]

Any such triangle is obtained exactly once.\\

Zimhoni proved this theorem skillfully by using the SternBrocot tree [3],[4].

Next we show the following theorem proved here [5]:\\

\noindent Theorem 2. 

A non-equilateral integer sided triangle is primitive eisensteinian if, and only if, it is of the class of triangles whose sides are ${\rm m^2{+}mn{+}n^2}$, ${\rm m^2{+}2mn}$ and either ${\rm n^2{+}2mn}$ or ${\rm m^2{-}n^2}$, where m and n are positive relatively prime integers with m>n>0 and m$\nequiv$n(mod 3). Any such triangle is achieved exactly once.\\

Zimhoni proved Theorem 1 by the aid of above-mentioned Theorem 2.\\

\noindent 2. For the case of  a $120^{\circ}$ angle\\

We originally define a definition as below: \\

\noindent Definition 2. 

A primitive sub-Eisensteinian triangle is an integer sided triangle with side lengths (a,b,c) such that the angle opposite the side of length a is of $120^{\circ}$ and gcd(a,b,c){=}1.\\

We assume a theorem and shall prove it as follows:\\

\noindent Theorem 3. 

There exists a one-to-one mapping between the set of non-equilateral primitive Eisensteinian triangles and the set of primitive sub-Eisensteinian triangles.

\begin{proof}
\quad\par

According to Defintion 1 and Theorem 1, a non-equilateral primitive Eisensteinian triangle is an integer triangle with side lengths (a,b,c) such that the angle opposite the side of length a is of $60^{\circ}$ and gcd(a,b,c){=}1 and b${>}$c.

Now we use a 3$\times$3 matrix
\[\rm
  S = \left(
    \begin{array}{ccc}
      1 & 0 & 0  \\
      0 & 0 & 1  \\
      0 & 1 & -1
    \end{array}
  \right) 
\]
 and present above-mentioned (a,b,c) as a column, then the product of the two is:
\[\rm
  S 
\left(
    \begin{array}{c}
    \rm a \\
    \rm  b \\
    \rm  c
    \end{array}
  \right) 
= \left(
    \begin{array}{ccc}
      1 & 0 & 0  \\
      0 & 0 & 1  \\
      0 & 1 & -1
    \end{array}
  \right) 
\left(
    \begin{array}{c}
    \rm  a \\
    \rm  b \\
    \rm  c
    \end{array}
  \right)
= 
\left(
    \begin{array}{c}
     \rm a \\
     \rm c \\
     \rm b-c
    \end{array}
  \right)
\]

Now we assume that (a,b,c) is an element the set of non-equilateral primitive Eisensteinian triangles.

By Definition 1, the Law of Cosines says that:\\

(2.1) \ \ \ \ \ \ ${\rm a^2=b^2+c^2-2bc\ cos60^\circ=b^2+c^2-bc} $\\

This is followed by:\\

(2.2) \ \ \ \ \ \ ${\rm a^2=c^2+(b-c)^2+c(b-c)=c^2+(b-c)^2-2c(b-c)\ cos120^\circ} $\\

According to Definition 2, (a,c,b${-}$c) is an element of the set of primitive sub-Eisensteinian triangles, because b${-}$c is a positive integer and gcd(a,c,b${-}$c)=1, owing to gcd(a,b,c)=1. \\

Inversely we assume that (a,b,c) is an element the set of primitive sub-Eisensteinian triangles.
Now we compute the inverse of S:
\[\rm
  S^{-1} = \left(
    \begin{array}{ccc}
      1 & 0 & 0  \\
      0 & 1 & 1  \\
      0 & 1 & 0
    \end{array}
  \right) 
\]\\
 and present above-mentioned (a,b,c) as a column, then the product of the two is:
\[\rm
  S^{-1} 
\left(
    \begin{array}{c}
     \rm a \\
     \rm b \\
     \rm c
    \end{array}
  \right) 
= \left(
    \begin{array}{ccc}
      1 & 0 & 0  \\
      0 & 1 & 1  \\
      0 & 1 & 0
    \end{array}
  \right) 
\left(
    \begin{array}{c}
     \rm a \\
     \rm b \\
     \rm c
    \end{array}
  \right)
= 
\left(
    \begin{array}{c}
     \rm a \\
     \rm b+c \\
     \rm b
    \end{array}
  \right)
\]

Now we assume that (a,b,c) is an element the set of primitive sub-Eisensteinian triangles.

By Definition 2, the Law of Cosines says that\\

(2.3) \ \ \ \ \ \ ${\rm a^2=b^2+c^2-2bc\ cos120^\circ=b^2+c^2+bc} $\\

This is followed by\\

(2.4) \ \ \ \ \ \ ${\rm a^2=(b+c)^2+b^2-(b+c)b=(b+c)^2+b^2-2(b+c)b\ cos60^\circ} $\\

According to Definition 1, (a,b${+}$c,b) is an element of the set of non-equilateral primitive Eisensteinian triangles , because b${+}$c${>}$b and gcd(a,b${+}$c,b)=1, owing to gcd(a,b,c)=1.\\

These complete the proof.

\end{proof}

By using Theorem 3, we can easily see as follows:\\

\noindent Corollary 1. 

An integer sided triangle is primitive sub-eisensteinian if, and only if, it is of the class of triangles whose sides are ${\rm m^2{+}mn{+}n^2}$, ${\rm n^2{+}2mn}$ and ${\rm m^2{-}n^2}$, where m and n are positive relatively prime integers with m>n>0 and m$\nequiv$n(mod 3). Any such triangle is achieved exactly once. 

\begin{proof}
\quad\par

By using S, we see
\[\rm
  S 
\left(
    \begin{array}{c}
    \rm m^2{+}mn{+}n^2 \\
    \rm m^2{+}2mn \\
    \rm n^2{+}2mn
    \end{array}
  \right) 
= \left(
    \begin{array}{ccc}
      1 & 0 & 0  \\
      0 & 0 & 1  \\
      0 & 1 & -1
    \end{array}
  \right) 
\left(
    \begin{array}{c}
    \rm m^2{+}mn{+}n^2 \\
    \rm m^2{+}2mn \\
    \rm n^2{+}2mn
    \end{array}
  \right)
= 
\left(
    \begin{array}{c}
     \rm m^2{+}mn{+}n^2 \\
     \rm n^2{+}2mn \\
     \rm m^2{-}n^2
    \end{array}
  \right)  \ and
\]
\[\rm
  S 
\left(
    \begin{array}{c}
    \rm m^2{+}mn{+}n^2 \\
    \rm m^2{+}2mn \\
    \rm m^2{-}n^2
    \end{array}
  \right) 
= \left(
    \begin{array}{ccc}
      1 & 0 & 0  \\
      0 & 0 & 1  \\
      0 & 1 & -1
    \end{array}
  \right) 
\left(
    \begin{array}{c}
    \rm m^2{+}mn{+}n^2 \\
    \rm m^2{+}2mn \\
    \rm m^2{-}n^2
    \end{array}
  \right)
= 
\left(
    \begin{array}{c}
     \rm m^2{+}mn{+}n^2 \\
     \rm m^2{-}n^2 \\
     \rm n^2{+}2mn
    \end{array}
  \right) .
\]

According to what Theorem 3 says, this completes Corollary 1.

\end{proof}

\noindent Theorem 4. 

An integer triangle is primitive sub-Eisensteinian if, and only if it has side lengths (a,b,c) or (a,c,b) such that:
\[\rm
  \left(
    \begin{array}{c}
    \rm a \\
    \rm b \\
    \rm c 
    \end{array}
  \right) = N\cdot v,\ where\ v{\in}
\left\{
\left(
    \begin{array}{c}
      7  \\
      5  \\
      3 
    \end{array}
  \right){,}
\left(
    \begin{array}{c}
      13  \\
      7  \\
      8 
    \end{array}
  \right)
  \right\}
\]

And N is a finite (or empty) product of the matrices:
\[\rm
  N_1=\left(
    \begin{array}{ccc}
      7 &  0 & -6 \\
      4 & -1 & -4 \\
      4 & 1 & -3
    \end{array}
  \right){,}
N_2=\left(
    \begin{array}{ccc}
      7 & 0 & -6 \\
      4 & -1 & 3 \\
      4 & 1 & 4
    \end{array}
  \right){,}
N_3=\left(
    \begin{array}{ccc}
      7 & 6 & 6 \\
      4 & 4 & 3 \\
      4 & 3 & 4
    \end{array}
  \right){,}
\]
\[\rm
  N_4=\left(
    \begin{array}{ccc}
      7 & 6 & 0 \\
      4 & 4 & 1 \\
      4 & 3 & -1
    \end{array}
  \right){,}
N_5=\left(
    \begin{array}{ccc}
      7 & -6 & 0 \\
      4 & -3 & 1 \\
      4 & -4 & -1
    \end{array}
  \right)
\]

Any such triangle is obtained exactly once.\\

\begin{proof}
\quad\par
According to Theorem 1, a non-equilateral integer sided triangle is primitive Eisensteinian if, and only if it has side lengths (a,b,c) or (a,b,b${-}$c) such that:
\[\rm
  \left(\rm
    \begin{array}{c}
    \rm a  \\
    \rm  b  \\
    \rm  c 
    \end{array}
  \right) = M\cdot v,\ where\ v{\in}
\left\{
\left(
    \begin{array}{c}
      7  \\
      8  \\
      5 
    \end{array}
  \right){,}
\left(
    \begin{array}{c}
      13  \\
      15  \\
       7 
    \end{array}
  \right)
  \right\}
\]\\
And M is a finite (or empty) product of the matrices:
\[\rm
  M_1=\left(
    \begin{array}{ccc}
      7 & -6 & 6 \\
      8 & -7 & 7 \\
      4 & -4 & 3
    \end{array}
  \right){,}
M_2=\left(
    \begin{array}{ccc}
      7 & 6 & -6 \\
      8 & 7 & -7 \\
      4 & 3 & -4
    \end{array}
  \right){,}
M_3=\left(
    \begin{array}{ccc}
      7 & 6 & 0 \\
      8 & 7 & 0 \\
      4 & 3 & 1
    \end{array}
  \right){,}
\]
\[\rm
  M_4=\left(
    \begin{array}{ccc}
      7 & 0 & 6 \\
      8 & 0 & 7 \\
      4 & 1 & 3
    \end{array}
  \right){,}
M_5=\left(
    \begin{array}{ccc}
      7 & 0 & -6 \\
      8 & 0 & -7 \\
      4 & 1 & -4
    \end{array}
  \right)
\]
\\
and any such triangle is obtained exactly once.

We now temporally consider the case of side lengths (a,b,c).

Then we set\\

(2.5) \ \ M=${\rm M_{p1}}$${\rm M_{p2}}$$\cdots$${\rm M_{pn}}$

\ \ \ \ \ \ \ \ \ \ \ \ where for any integer k(1$\leq$k$\leq$n) and

\ \ \ \ \ \ \ \ \ \ \ \ \ \ \ ${\rm M_{pk}}$ = ${\rm M_1}$ or ${\rm M_2}$ or ${\rm M_3}$ or ${\rm M_4}$ or ${\rm M_5}$.\\

Moreover we obtain that:\\

(2.6) \ \ S$\cdot$M$\cdot$v=S(${\rm M_{p1}}$${\rm M_{p2}}$$\cdots$${\rm M_{pn}}$)v

\ \ \ \ \ \ \ =(${\rm SM_{p1}S^{-1}}$)(${\rm SM_{p2}S^{-1}}$)$\cdots$(${\rm SM_{pn}S^{-1}}$)(S$\cdot$v)\\

According to the proof of Theorem 2, S$\cdot$v is an element of the set of primitive sub-Eisensteinian triangles. And we easily see that two different elements of the set of primitive Eisensteinian triangles is respectively mapped to different elements of the set of primitive Eisensteinian triangles by S. 

And owing to Theorem 1, if any two certificates for M that have a different expression from each other, then their values are also different from each other. 

We arrange (2.6) as follows:\\

(2.7) \ \ S$\cdot$M$\cdot$v=S(${\rm M_{p1}}$${\rm M_{p2}}$$\cdots$${\rm M_{pn}}$)v

\ \ \ \ \ \ \ \ \ \ =(${\rm SM_{p1}S^{-1}}$)(${\rm SM_{p2}S^{-1}}$)$\cdots$(${\rm SM_{pn}S^{-1}}$)(S$\cdot$v)

\ \ \ \ \ \ \ \ \ \ =(${\rm N_{p1}}$)(${\rm N_{p2}}$)$\cdots$(${\rm N_{pn}}$)(S$\cdot$v)

\ \ \ \ \ \ \ \ \ \ =${\rm N}$(S$\cdot$v).\\

We let:\\

(2.8) \ \ N=${\rm N_{p1}}$${\rm N_{p2}}$$\cdots$${\rm N_{pn}}$

\ \ \ \ \ \ \ \ \ \ \ \ where for any integer k(1$\leq$k$\leq$n) and

\ \ \ \ \ \ \ \ \ \ \ \ \ \ \ ${\rm N_{pk}}$ = ${\rm SM_{1}S^{-1}}$ or ${\rm SM_{2}S^{-1}}$ or ${\rm SM_{3}S^{-1}}$ or ${\rm SM_{4}S^{-1}}$ or ${\rm SM_{5}S^{-1}}$.\\

In the above-mentioned equations, any two certificates for N that have a different expression from each other, then their values are also different from each other.

Then by computing we see: 

\[\rm
  N_1=\left(
    \begin{array}{ccc}
      1 & 0 & 0 \\
      0 & 0 & 1 \\
      0 & 1 & -1
    \end{array}
\right)
\left(
    \begin{array}{ccc}
      7 & -6 & 6 \\
      8 & -7 & 7 \\
      4 & -4 & 3
    \end{array}
  \right)
\left(
    \begin{array}{ccc}
      1 & 0 & 0 \\
      0 & 1 & 1 \\
      0 & 1 & 0
    \end{array}
  \right){=}
\left(
    \begin{array}{ccc}
      7 & 0 & -6 \\
      4 & -1 & -4 \\
      4 & 1 & -3
    \end{array}
  \right){,}
\]
\[\rm
  N_2=\left(
    \begin{array}{ccc}
      1 & 0 & 0 \\
      0 & 0 & 1 \\
      0 & 1 & -1
    \end{array}
\right)
\left(
    \begin{array}{ccc}
      7 & 6 & -6 \\
      8 & 7 & -7 \\
      4 & 3 & -4
    \end{array}
  \right)
\left(
    \begin{array}{ccc}
      1 & 0 & 0 \\
      0 & 1 & 1 \\
      0 & 1 & 0
    \end{array}
  \right){=}
\left(
    \begin{array}{ccc}
      7 & 0 & 6 \\
      4 & -1 & 3 \\
      4 & 1 & 4
    \end{array}
  \right){,}
\]
\[\rm
  N_3=\left(
    \begin{array}{ccc}
      1 & 0 & 0 \\
      0 & 0 & 1 \\
      0 & 1 & -1
    \end{array}
\right)
\left(
    \begin{array}{ccc}
      7 & 6 & 0 \\
      8 & 7 & 0 \\
      4 & 3 & 1
    \end{array}
  \right)
\left(
    \begin{array}{ccc}
      1 & 0 & 0 \\
      0 & 1 & 1 \\
      0 & 1 & 0
    \end{array}
  \right){=}
\left(
    \begin{array}{ccc}
      7 & 6 & 6 \\
      4 & 4 & 3 \\
      4 & 3 & 4
    \end{array}
  \right){,}
\]
\[\rm
  N_4=\left(
    \begin{array}{ccc}
      1 & 0 & 0 \\
      0 & 0 & 1 \\
      0 & 1 & -1
    \end{array}
\right)
\left(
    \begin{array}{ccc}
      7 & 0 & -6 \\
      8 & 0 & -7 \\
      4 & 1 & 3
    \end{array}
  \right)
\left(
    \begin{array}{ccc}
      1 & 0 & 0 \\
      0 & 1 & 1 \\
      0 & 1 & 0
    \end{array}
  \right){=}
\left(
    \begin{array}{ccc}
      7 & 6 & 0 \\
      4 & 4 & 1 \\
      4 & 3 & -1
    \end{array}
  \right){,}
\]
\[\rm
  N_5=\left(
    \begin{array}{ccc}
      1 & 0 & 0 \\
      0 & 0 & 1 \\
      0 & 1 & -1
    \end{array}
\right)
\left(
    \begin{array}{ccc}
      7 & 0 & -6 \\
      8 & 0 & -7 \\
      4 & 1 & -4
    \end{array}
  \right)
\left(
    \begin{array}{ccc}
      1 & 0 & 0 \\
      0 & 1 & 1 \\
      0 & 1 & 0
    \end{array}
  \right){=}
\left(
    \begin{array}{ccc}
      7 & -6 & 0 \\
      4 & -3 & 1 \\
      4 & -4 & -1
    \end{array}
  \right),
\]
\[\rm Sv=
\left(
    \begin{array}{ccc}
      1 & 0 & 0 \\
      0 & 0 & 1 \\
      0 & 1 & -1
    \end{array}
  \right)
\left(
    \begin{array}{c}
      7  \\
      8  \\
      5 
    \end{array}
  \right)
{=}
\left(
    \begin{array}{c}
      7  \\
      5  \\
      3 
    \end{array}
  \right) \ \ or
\]
\[\rm \ \ \ =
\left(
    \begin{array}{ccc}
      1 & 0 & 0 \\
      0 & 0 & 1 \\
      0 & 1 & -1
    \end{array}
  \right)
\left(
    \begin{array}{c}
      13  \\
      15  \\
      7 
    \end{array}
  \right)
{=}
\left(
    \begin{array}{c}
      13  \\
      7  \\
      8 
    \end{array}
  \right).
\]

The general expression of the elements of the set of primitive sub-Eisensteinian triangles, mapped by S from the elements of the set of non-equilateral primitive Eisensteinian triangles, in the case of side lengths (a,b,c), is (2.7).

Then next we consider the case of side lengths (a,b,b${-}$c). In this connection, the pair (a,b,c) and (a,b,b${-}$c) is called the the twin Eisensteinian triplet [1].

With the aid of Theorem 1, we obtain:
\[\rm
(2.9)\ \ S\left(
    \begin{array}{c}
      \rm a \\
      \rm b \\
      \rm b-c
    \end{array}
  \right)=
S\left(
    \begin{array}{ccc}
      1 & 0 & 0 \\
      0 & 1 & 0 \\
      0 & 1 & -1
    \end{array}
  \right)
\left(
    \begin{array}{c}
      \rm a \\
      \rm b \\
      \rm c
    \end{array}
  \right)=
S\left(
    \begin{array}{ccc}
      1 & 0 & 0 \\
      0 & 1 & 0 \\
      0 & 1 & -1
    \end{array}
  \right)(M\cdot v)
\]
\[\rm
=S\left(
\begin{array}{ccc}
      1 & 0 & 0 \\
      0 & 1 & 0 \\
      0 & 1 & -1
    \end{array}
  \right)S^{-1}(SMv)=
\left(
    \begin{array}{ccc}
      1 & 0 & 0 \\
      0 & 0 & 1 \\
      0 & 1 & 0
    \end{array}
  \right)(SMv)
\]

It is evident that: 
\[\rm 
(2.10) \ \ \ \ \ \ \   
\left(
    \begin{array}{ccc}
      1 & 0 & 0 \\
      0 & 0 & 1 \\
      0 & 1 & 0
    \end{array}
  \right)
\left(
    \begin{array}{c}
      \rm a \\
      \rm b \\
      \rm c
    \end{array}
  \right)=
\left(
    \begin{array}{c}
      \rm a \\
      \rm c \\
      \rm b
    \end{array}
  \right),\ \ \ \ \ \ \ \ \ \ \ \ \ \ \ \ \ \ \ \ \ \ \ \ 
\]
therefore we see as follows:

If the case of side lengths (a,b,c) for primitive Eisensteinian corresponds to the case of side lengths (a,b,c) for primitive sub-Eisensteinian, then the case of side lengths (a,b,b${-}$c) for primitive Eisensteinian corresponds to the case of side lengths (a,c,b) for primitive sub-Eisensteinian. 

Thus any primitive sub-Eisensteinian triangle is also obtained by pre-multiplication (7,5,3) or (13,7,8) by a product of five matrices.

These complete the proof.

\end{proof}
%
%
%

\centerline{References}
\noindent [1] Noam Zimhoni, arXiv, available at https://arxiv.org/abs/1904.11782 \\
\noindent [2] B Berggren. Pytagoreiska trianglar. Tidskrift for elementar matematik, fysik och kemi, 17:129-139, 1934.  \\
\noindent [3] Moritz A Stern. Ueber eine zahlentheoretische funktion. Journal fur die reine und angewandte Mathematik, 55:193-220, 1858. \\
\noindent [4] Achille Brocot. Calcul des rouages par approximation, nouvelle methode. Revue chronometrique, 3:186-194, 1861. \\
\noindent [5] Shin-ichi Katayama. Modified farey trees and pythagorean triples. Journal of mathematics, the University of Tokushima, 47, 2013.
\end{document}